\documentclass[11pt,twoside]{article}
\usepackage{graphicx} 
\usepackage{ulem}
\usepackage{amssymb}
\usepackage[mathscr]{eucal}
\usepackage{eufrak}
\usepackage{amsmath,amsthm}
\usepackage{physics}
\usepackage{mathrsfs}
\usepackage{lgreek}
\usepackage[shortlabels]{enumitem}
\usepackage[colorlinks=true,
   urlcolor=blue,           filecolor=green,      
   citecolor=green,      
   linkcolor=red,           bookmarks=true,
  unicode,
   plainpages=false,   ]{hyperref}
\usepackage{color}
\font\bg=cmbx10 scaled 1300

\definecolor{royalblue}{rgb}{0,0,0.128}
\def\bp{\begin{proof}}
\def\ep{\end{proof}}
\def\n{\nabla}

\def\sfrac#1#2{\mbox{\Large$\frac{#1}{#2}$}}
\def\intl#1{\int\limits_{#1}}
\def\intll#1#2{\int\limits_{#1}^{#2}}
\def\dm{|\hskip-0.05cm|}

\def\OO{\Omega}
\def\displ{\displaystyle}
\def\VSE{\vspace{6pt}\\&\displ }
\def\VS{\vspace{6pt}\\\displ }
\def\rf#1{{\rm(\ref{#1})}}

\def\R{\Bbb R}

\def\à{à}

\def\vep{\varepsilon}

\def\be{\begin{equation}}
\def\ba{\begin{array}}
\def\ea{\end{array}}
\def\ee{\end{equation}}
\def\vs1{\vspace{1ex}}

\def\po{\partial\Omega}

\def\Ã©{\'{e}}
\def\Ãš{\`{e}}
\def\rot{{\rm curl}\hskip0.02cm}
\newtheorem{lemma}
{\bf Lemma} 
\font\sc=cmcsc10
\date{\today}

\newtheorem{defi}
{\bf Definition} 
\newtheorem{tho}
{\bf Theorem} 
\newtheorem{rem}
{\sc Remark} 
 
\newtheorem{prop}
{\bf Proposition} 

\title{\bg The motion of a rigid body in a viscous fluid:\\new results for strong solutions, uniqueness and integrability properties}
\author{\sc Paolo Maremonti and Filippo Palma
\thanks{Dipartimento di Matematica e Fisica,  
Universit\`{a} degli 
Studi della Campania
``L. Vanvitelli'', via Vivaldi 43, 81100 \null\hskip0.55cmCaserta,
 Italy.\newline\null\hskip0.55cm
paolo.maremonti@unicampania.it
\newline\null\hskip0.55cmfilippo.palma@unicampania.it
\newline\null\hskip0.55cm The  research activity  is performed under the
auspices of   GNFM-INdAM. }}
\begin{document}
\maketitle
\begin{abstract}
    In this note, we show two results in the setting of Galdi-Silvestre strong solutions for the rigid body-viscous fluid interaction. The former, under an additional integrability assumption on the gradient of the initial datum, proves that the time derivative of the solution belongs to $L^2(0,T;L^2(\Omega))$. The latter, thanks to a further assumption only on one solution, proves that the uniqueness holds in the quoted setting. However, our extra assumption for the uniqueness is certainly verified under the integrability assumption on the gradient of the initial datum. Hence, the set of solutions enjoying the uniqueness is not empty.
\end{abstract}
\section{Introduction}
Let us consider the motion of a rigid body $\mathcal{B}$ in a viscous fluid $\mathcal{F}$ that fills the entire three-dimensional space exterior to $\mathcal{B}$. The shape of $\mathcal{B}$ can be arbitrary, we only require the boundary $\partial\mathcal{B}$ to be sufficiently smooth. We investigate  the motion of the system $\mathbf{S}:=\{\mathcal{B},\mathcal{F}\}$.\par In a reference frame attached to the body $\mathcal{B}$, with origin at its center of mass $G$, we denote by $\Omega := \mathbb{R}^3 \setminus \mathcal{B}$ the region occupied by the fluid. In the absence of external forces and torques, the equations governing the motion of the fluid within $\Omega$ are given by:
\begin{equation} \label{eq:model}
    \begin{cases}
        & u_t-\Delta u= - [(u-V)\cdot \nabla u+\omega\times u]-\nabla\pi\quad \forall (t,x) \in (0,T)\times \Omega\,, \\
        & \nabla \cdot u=0 \quad \forall (t,x) \in (0,T)\times \Omega\,, \\
        & u(t,x)=V(t,x)=\xi(t)+\omega(t)\times x \quad \forall (t,x) \in (0,T)\times \partial \Omega\,, \\
        & \lim_{\abs{x}\to \infty} u(t,x)=0 \quad \forall t\in (0,T)\,, \\
        &\dot{\xi} + \omega \times \xi +  \int_{\partial\Omega} \mathbb{T}(u, \pi)\cdot \nu=0 \quad \forall t\in (0,T)\,,\\
        & I \cdot \dot{\omega}+ \omega \times(I\cdot \omega)+ \int_{\partial\Omega} x \times \mathbb{T}(u,\pi) \cdot \nu=0 \quad \forall t\in (0,T)\,,\\
        & \xi(0)= \xi_0\,, \quad \omega(0)=\omega_{0}\,, \\
        & u(0,x)=u_0(x)\,, \quad \forall x \in \Omega\,.
    \end{cases}
\end{equation}
Here, $u:(0,T)\times \Omega \to \mathbb{R}^3$ is the velocity of the fluid, $\pi:(0,T)\times \Omega \to \mathbb{R}$ is the pressure field, $V:(0,T)\times \bar{\OO}\to \mathbb{R}^3$ is the velocity  of the rigid body expressed by $V(t,x)=\xi(t)+\omega(t)\times x$, where  $\xi:(0,T)\to\mathbb{R}^3$ is the translation velocity and $\omega:(0,T)\to\mathbb{R}^3$ is the rotation velocity, $I$ is the inertia tensor, $\nu$ is the outer unit normal vector of the surface $\partial\Omega$ and $\mathbb{T}$ is the Cauchy stress tensor of a Newtonian viscous fluid, defined by
\[
\mathbb{T}(u,\pi):= -\pi\mathbb{I}+2\mu\mathbb{D}(u),
\]
where $\mathbb{I}$ is the identity tensor, $\mu$ is the  dynamic viscosity and $\mathbb{D}(u)=\frac{1}{2}(\nabla u+ (\nabla u)^T)$ is the symmetric part of $\nabla u$. For the sake of simplicity, we assume $\mu=1$. For an exhaustive overview of the model, we refer to the seminal article by Galdi \cite{Ga:Hand}. \par
The model \eqref{eq:model} has been extensively investigated in the last decades by several authors. \par First results concerning the existence of global weak solutions go back to Serre \cite{Serre}. \par A first result on strong solutions was obtained by Galdi and Silvestre \cite{GS}, in their work it was shown that given a sufficiently smooth $u_0$ with the compatibility condition of rigid motion on the boundary $\partial\Omega$, there exists at least one local strong solution to problem \eqref{eq:model}. The proof, based on the Faedo Galerkin method and the invading domains technique, was given in the $L^2$-setting. \par In the $L^p$-setting, with the assumption of a purely rotating body ($\xi\equiv 0$), the model was investigated by Geissert, Heck and Hieber \cite{Gei-Hec-Hie}, who showed the existence of a unique local mild solution employing Kato's iteration procedure. \par Recently, Galdi \cite{Ga:23} proved that, if the initial datum satisfies a suitable smallness condition, there exists at least one strong global solution which tends to the rest state asymptotically in time. \par  An improvement concerning the asymptotic behavior in time of a global solution was obtained by Galdi and Maremonti \cite{Gal-Mar}.  \par In both the papers \cite{GS} and \cite{Ga:23}, the uniqueness of the strong solution is still an open problem.\par Concerning the steady-state problem associated with system $\mathbf{S}$, we refer to \cite{GS-sta}. We also quote the work by Hishida \cite{His} on the steady motion of a rigid body around a rotating obstacle. \par In a different frame, that allows to ``eliminate'' the critical term $\omega\times x$ in the equation \eqref{eq:model}$_1$, a $L^p$-theory for the motion of a rigid body in a Newtonian or a generalized Newtonian viscous fluid was developed by Geissert, G\"{o}tze and Hieber \cite{Gei-Got-Hie}. They proved a result of strong local well-posedness using a fixed-point argument. \par In the same frame, a first result on global existence and uniqueness of a strong solution for ``small'' initial data was proved by Cumsille and Takahashi \cite{Cum-Tak}, while large time behavior of the solution was investigated by Ervedoza et al. \cite{Erv-Mai-Tuc}.\par   The purpose of this note is twofold. On one hand we find a sufficient condition on the initial datum in order to get $u_t\in L^2(0,T;L^2(\Omega))$. On the other hand, we find a sufficient condition so that a strong solution in the setting of Galdi-Silvestre is unique. However, we stress that the sufficient condition for the uniqueness is certainly satisfied as soon as we consider the same assumption on the initial datum that allows us to obtain the property $u_t\in L^2(0,T;L^2(\OO))$. \par
Before stating our main results, we briefly introduce the notation for the function spaces which will be used throughout the paper. We start by recalling that by $L^p(\Omega)$ and $W^{m,p}(\Omega)$ ($m\in\mathbb{N}$, $p\ge1$) we refer to the usual Lebesgue and Sobolev spaces, whose norms are denoted by $\norm{\cdot}_p$ and $\norm{\cdot}_{m,p}$, respectively. Furthermore, we denote by $L^p(\Omega, \abs{x}^{a})$, $a\in \mathbb{N}_0$, the following weighted Lebesgue space
\[
L^p(\Omega,\abs{x}^a):=\{u \, : \, \abs{x}^a u\in L^p(\Omega)\}.
\]
Setting
\[
\mathscr{C}_0(\Omega):=\{u \in C_0^{\infty}(\Omega) \, : \nabla \cdot u=0\},
\]
we define the space $J^p(\Omega)$, $p\in (1,+\infty)$, as the closure of $\mathscr{C}_0(\Omega)$ with respect to the norm $\norm{\cdot}_p$\,. \par
Let $X$ be a Banach space. The spaces $C([a,b];X)$ and $L^{p}(a,b;X)$ are defined as follows
\[
C([a,b];X):=\{u:[a,b]\to X \, : \, u\text{ is continuous in the norm} \, \, \norm{\cdot}_X\};
\]
\[
L^{p}(a,b;X):=\{ u:[a,b]\to X \, :\,  \biggl(\int_a^b \norm{u(t)}_{X}^p \, dt\biggr)^{\frac{1}{p}}<\infty \}. 
\]
Since we need to deal with the motion of both a rigid body and a viscous fluid, we should introduce some function spaces that allow us to follow the evolution of both elements of the system $\mathbf{S}$.\par The symbol $\mathbb E_3$ denotes the three-dimensional vector space on the field $\R$. We set
\[
\mathcal{R}:=\{\bar{u} \in C^{\infty}(\mathbb{R}^3)\, : \, \bar{u}(x)=\bar{u}_1+\bar{u}_2\times x, \,\, \bar{u}_1,\bar{u}_2\in \mathbb{E}_3\}.
\]
The vectors $\bar{u}_1$ and $\bar{u}_2$ are called \textit{characteristic vectors} of the rigid motion $\bar{u}$. \par The vector space $\mathcal{H}(\Omega)$ and its subspace $\mathcal{V}(\Omega)$ are defined as
\be\label{HSP}
\mathcal{H}(\Omega):= \{u \in L^2(\Omega) \,:\, \nabla \cdot u=0, \, (u-\bar{u})\cdot \nu|_{\partial\Omega}=0, \, \, \text{ for some }\bar{u}\in \mathcal{R}\}
\ee
and
\be\label{VSP}
\mathcal{V}(\Omega):=\{u\in W^{1,2}(\Omega) \, : \, \nabla \cdot u=0, \, u|_{\partial\Omega}=\bar{u}\, \, \text{ for some }\bar{u}\in \mathcal{R}\}.
\ee
We introduce the scalar products
\be\label{SPH}
(u,v)_{\mathcal{H}(\Omega)}:= \bar{u}_1\cdot \bar{v}_1+\bar{u}_2\cdot \bar{v}_2+ \int_{\Omega}(u,v) \, dx
\ee
and
\be\label{SPV}
(u,v)_{\mathcal{V}(\Omega)}:= \bar{u}_1\cdot \bar{v}_1+\bar{u}_2\cdot \bar{v}_2+ \int_{\Omega}(u,v) \, dx+\int_{\Omega}(\nabla u,\nabla v) \, dx,
\ee
where $(\cdot,\cdot)$ is the standard $L^2$ scalar product, respectively in $\mathcal H(\OO)$ and $\mathcal V(\OO)$.\par 
The two spaces, endowed with their scalar products, are Hilbert spaces. \par For more details about function spaces, we refer to \cite{Adams}, for the functional setting concerning the model, we quote \cite{Ga:Hand, GS}. \par
 We remark that, throughout the following, the symbols $c$ and $C$ will be used to denote generic non-negative constants, whose values may vary from line to line.\par
Furthermore, we denote by $\OO_R$ the open subset $\OO \cap \{x \in \mathbb{R}^3 : |x| < R\}$, and by $\OO^R$ the open subset $\OO \cap \{x \in \mathbb{R}^3 : |x| > R\}$.
\begin{defi}\label{GSS}{\sl We say that the system $\mathbf S$ has a regular motion  in the sense of Galdi-Silvestre, if we have that the fields $(u,\pi,V)$   solve system \rf{eq:model} a.e. on  $(0,T)\times\OO$ and they enjoy of the following properties:\begin{equation}\label{GSS-I}
\begin{aligned}
&u\in L^{\infty}(0,T; \mathcal V(\Omega)) \cap L^2(0,T; W^{2,2}(\Omega))\,, \\ &\nabla\pi \in L^2(0,T; L^2(\Omega))\,,\\
& \xi,\omega \in W^{1,2}(0,T)\,, \\
& u_t \in L^2(0,T; L^2(\Omega_R))\,, \quad \forall R>\text{diam}(\mathcal{B})\,.
\end{aligned}
\end{equation}
Moreover, it holds
\begin{equation}\label{GSS-II}
    \begin{aligned}
       & \xi\,,\omega \in C([0,T'])\,, \,\, \forall \,T'<T\,, \,\, \xi(0)=\xi_0\,, \,\, \omega(0)=\omega_{0}\,, \\
       & u\in C([0,T']; W^{1,2}(\Omega_R))\,, \,\, \forall\, T'<T\,\, \text{and } \forall \,R>\text{diam}(\mathcal{B})\,, \,\, u(0,\cdot)=u_0(\cdot)\,.
    \end{aligned}\ee
}\end{defi}
The previous definition is consistent by virtue of the following result:
\begin{tho} \label{thm:esistenza}Let $(\xi_0,\omega_0)\in \mathbb E_3\times\mathbb E_3$.
    Let $u_0 \in \mathcal{V}(\Omega)$ be such that $u_0|_{\partial\Omega}= \xi_0+ \omega_0\times x$\,. Then, there exists at least one solution in the sense of Definition\,\ref{GSS}. Additionally, there exists $\delta>0$ such that, if
    \begin{equation}
        \norm{u_0}_{1,2}+ \abs{\xi_0}+\abs{\omega_0} \le \delta\,,
    \end{equation}
   then $T=\infty$, and the following property of decay holds:
\begin{equation}
    \lim_{t\to +\infty} (\norm{u(t)}_6+ \norm{\nabla u(t)}_2+\abs{\xi(t)}+ \abs{\omega(t)})=0\,.
\end{equation}
\end{tho}
The reader can check the proof of the previous theorem in \cite{GS} for the local existence and in \cite{Ga:23} for the global existence. \begin{rem}{\rm We want to stress that the results by Galdi-Silvestre and Galdi are obtained considering that the body could be subject to prescribed external forces $F$ and torques $M$. Since the presence of prescribed external forces and torques does not add analytical issues related to our questions, for the sake of brevity we set $F=M=0$.}\end{rem} \par
\begin{rem}
    {\rm In \cite{GS} it was shown that $\nabla\pi \in L^2(0,T; L^2(\Omega_R))$ for all $R>\text{diam}(\mathcal{B})$, the improved property $\nabla\pi \in L^2(0,T; L^2(\Omega))$ was first proved for a time-periodic solution in \cite{GS-per}, however, it can easily be extended to the local solution deduced in \cite{GS}}.
\end{rem}
We are now in position to state our main results.
\begin{tho}\label{TDU}{\sl Let $(u_0,\xi_0,\omega_0)\in \mathcal V(\OO)\times \mathbb{E}_3\times\mathbb{E}_3$, with ${u_0}_{|{\po}}=\xi_0+\omega_0\times x$, and $\n u_0\in L^2(\Omega,\abs{x})$. Then there exists a strong solution to problem\,\eqref{eq:model} enjoying properties\,\eqref{GSS-I}$_{1,3,4}$, \eqref{GSS-II} and  
\be \label{GSS-old}
\n \pi \in L^2(0,T;L^2(\Omega_R))\,, \quad \forall R>\text{diam}(\mathcal{B})\,,
\ee
such that 
\be\label{TDU-I}\hskip0.05cm\n u(t, \cdot)\in L^2(\Omega,\abs{x})\,,\mbox{ for all }t\in[0,T)\,,\mbox{\; and \;}\intll0T\dm |x|D^2u\dm_2^2dt<\infty\,. \ee Moreover, if $(u,\pi,V)$ satisfies also \eqref{GSS-I}$_2$, then
\be \label{TDU-II} u_t\in L^2(0,T;L^2(\OO))\,.\ee}\end{tho}
\begin{rem}{\rm Employing just \rf{GSS-old},
     Theorem\,\ref{TDU}, in particular, suggests that the property \eqref{GSS-I}$_2$ for the pressure field and the weighted integrability property for the gradient in \eqref{TDU-I} are completely independent. \par Differently, showing that $u_t\in L^2(0,T;L^2(\Omega))$ without assuming \eqref{GSS-I}$_2$ seems still out of reach.}
\end{rem}
\begin{rem}
   {\rm Theorem\,\ref{TDU} allows us to state that a strong solution $(u,\pi,V)$, ensured by Theorem\,\ref{GSS}, corresponding to an initial datum $u_0\in\mathcal{V}(\Omega)$ enjoying the additional property $\n u_0\in L^2(\Omega,\abs{x})$, fulfills the following requirement:
    \[
    u\in C([0,T'];W^{1,2}(\Omega))\,,
    \]
    for all $T'<T$\,.}
\end{rem}
\begin{tho} \label{UT}{\sl 
    Let $(u,\pi,V)$ be a solution to \eqref{eq:model} in the sense of Theorem\,\ref{thm:esistenza} corresponding to the initial datum $(u_0,\xi_0,\omega_0)\in\mathcal{V}(\Omega)\times \mathbb{E}_3\times\mathbb{E}_3$. Assume that 
    \be\label{UT-I}
    \hskip0.02cm\n u\in L^2(0,T; L^2(\Omega,\abs{x}))\,. \,\, 
    \ee
Then,  the solution is unique in the set of solutions of Definition\,\ref{GSS}.}
\end{tho}
\begin{rem}
    {\rm Theorem\,\ref{TDU} proves that condition \rf{UT-I} is satisfied at least for an initial datum $u_0\in \mathcal{V}(\Omega)$ such that $\n u_0\in L^2(\OO,|x|)$. We would like to stress that Theorem\,\ref{UT} can be interpreted in the sense of Leray-Serrin: one solution must be regular and satisfy \rf{UT-I}, while the other only needs to have sufficient regularity to ensure the validity of the energy relation for the equations of the difference between the two motions (see the proof of the theorem). For the sake of brevity, we omit the details of the proof of this sentence.}\end{rem}
\begin{rem}{\rm 
    Let us point out that estimate \eqref{TDU-II} suggests the following equivalence:
    \begin{equation}\label{eq:deduzione}
        u_t\in L^2(0,T;L^2(\Omega)) \iff \omega\times x\cdot \n u \in L^2(0,T;L^2(\Omega))\,.
    \end{equation}
    Roughly speaking, the equivalence \eqref{eq:deduzione} suggests that, a.e. in \( t > 0 \), the equations jointly with the
summability of \( u_t \) and with the regularity of \((u,\pi)\) imply a.e. in $t\in(0,T)$ a sort of stronger summability of \(\nabla u\),
that is the summability  in \( L^2(\Omega) \) with a suitable increasing weight, with respect to the simple
property \(\nabla u \in L^2(\Omega)\).
\par
Hence, \rf{eq:deduzione} implies
 \begin{equation}\label{eq:deduzione 1}
        u_t(t,\cdot)\in L^2(\Omega) \iff \omega(t)\times x\cdot\n u(t,\cdot) \in L^2(\Omega), \quad \text{a.e.  in }t\in(0,T).
    \end{equation}
\par 
If we assume a smooth initial datum in such a way that we can guess the validity of the
equivalence \eqref{eq:deduzione 1} up to the initial instant, then this means that some compatibility between the condition
and the smooth initial datum is due.
\par Since from a physical point of view the model of the fluid is independent of the time instant
of the motion that one considers, our assumption on the initial datum, although slightly stronger
than \(\omega(t) \times x \cdot \nabla u(t) \in L^2(\Omega)\) (moreover, in our setting only a sufficient condition for the result), appears
in some sense natural in the dynamic of system $\mathbf{S}$}.
\end{rem}
\section{Preliminary results}
Let $\text{diam}(\mathcal{B})<R_0\le R$. We denote by $\zeta_0 (x)$ and $ \zeta_1(x)$ two ordinary smooth cut-off functions such that $\zeta_0(x)=1$ for $|x|\leq R_0$ and $\zeta_1(x)=1$ for $|x|\leq R$, $\zeta_0(x)=0$ for $|x|\geq 2R_0$ and $\zeta_1(x)=0$ for $|x|\geq 2R$, and  their derivatives have support contained in the shell $R_0\leq |x|\leq 2R_0$   and in the shell $R\leq |x|\leq 2R$, respectively. \par We set \be\label{eq:zeta} \zeta(x):=(1- \zeta_0(x))\zeta_1(x), \mbox{ for all }x\in\R^3\,.\ee
Taking into account the definition of $\zeta_0$ and $\zeta_1$\,, we easily deduce that $\zeta(x)$ is smooth and $\n \zeta(x)$ and $\n\n \zeta(x)$ have compact support contained in the union of the shells $\{x:R\leq |x|\leq 2R\}\cup\{x:R_0\leq |x|\leq 2R_0\}$ with
$$|\n\zeta(x)|+R|\n\n\zeta(x)|\leq cR^{-1}\,,\,\mbox{ for all }|x|>R\,.$$ 
For a function $g\in C^3(\R^3)$ let us consider the identity
\be\label{FFL}\Delta g=\n\n\cdot g-\rot\rot g\,,\mbox{ for all }x\in\R^3\,,\ee
that also furnishes 
$$\Delta\rot g=-\rot\rot\rot g\,,\mbox{ for all }x\in\R^3\,.$$ Assuming that $g$ and its derivates belong to $L^2(\R^3)$, we get the identity
\be\label{IDM}\dm \n\rot g\dm_2=\dm \rot\rot g\dm_2\,.\ee
We set $g(x):=h(x)u(x)$\,, with $h$ scalar smooth function with compact support. The following identities hold:
$$\ba{l}\rot(hu)=h\rot u+\n h\times u\,,\VS \n\rot(h u)=h\n\rot u+\n h\otimes \rot u+[\vep_{ijk}D_{x_\ell}(u_iD_{x_j}h)]e_\ell\otimes e_k\,,\,i,j,k,\ell\in\{1,2,3\}\,,\VS
\rot\rot(hu)= h\rot\rot u+\n h\times\rot u+ \n h\n\cdot u-u\Delta h-\n h\cdot\n u+u\cdot\n\n h\,,
\ea$$ where $\vep_{ijk}$ are the Levi-Civita symbols and the Einstein sum convention is employed. Substituting the previous formulas in \rf{IDM} and handling suitably the term of the sum, one gets
\be\label{EDM}\dm h\n b\dm_2\leq \dm h\hskip0.02cm\rot b\dm_2+c\Big[\dm (|\n h|+|\n\n h|)u\dm_2+\dm |\n h|\n u\dm_2\Big]\,,\ee
where $b:=\rot u$ and $c$ is a non negative constant. \par We recall other inequalities of this kind.\par  Assume  $g\in C^3_0(\R^3)$\,. Assume  $|x|\rot g\in L^2(\R^3)$ and $|x|\n \cdot g\in L^2(\mathbb{R}^3)$. Then, we get
\be\label{RGW}\dm |x|\n g\dm_2\leq c\Big[\dm |x|\rot g\dm_2+\dm |x|\n\cdot g\dm_2\Big]\,,\ee with $c$ independent of $g$. Assume $|x|\n\rot g\in L^2(\OO)$ and $|x|\n\n\cdot g\in L^2(\OO)$, then we get \be\label{RDSW}\dm |x|D^2g\dm_2 \leq c\Big[\dm |x|\n\rot g \dm_2+\dm |x|\n \n\cdot g\dm_2\Big]\,,\ee with $c$ independent of $g$. \par We start again from \rf{FFL}. The function $g$ admits the representation formula $$g(x)=-\intl{\R^3}\n_y\mathcal E(x-y)\times\rot g(y)dy-\intl{\R^3}\n_y\mathcal E(x-y)\n\cdot g(y)dy\,.$$ Considering its gradient, we get 
$$ \n g(x)=-\n\intl{\R^3}\n_y\mathcal E(x-y)\times\rot g(y)dy-\n\intl{\R^3}\mathcal \n_y E(x-y)\n\cdot g(y)dy\,.$$ Hence, estimate \rf{EDM} is a consequence of the Calderon-Zygmund theorem with weights, see \cite{St}.  \par Analogously, applying two times the gradient operator to function $g$, we get
$$\n\n g(x)=\n\n\intl{R^3}\mathcal E(x-y)\rot\rot g(y)dy+\n\n\intl{\R^3}\mathcal E(x-y)\n\n\cdot g(y)dy\,.$$ Employing again the Calderon-Zygmund theorem with weights, and employing the pointwise estimate $|\rot\rot g(x)|\leq c|\n\rot g(x)|$\,, we arrive at \rf{RDSW} with $c$  independent of $g$\,.
\par 
As it is known, estimate \rf{EDM}  can be extended to the set of function $u\in W^{2,2}(\R^3)$. \par Furthermore,  estimate  \rf{RGW} to functions $u$ belonging to completion of $C_0^\infty(\R^3)$ with respect to the metric of  $|x|\rot u\in L^2(\R^3)$ and $|x|\n\cdot u\in L^2(\mathbb{R}^3)$. \par Finally, estimate \rf{RDSW} to functions $u$ belonging to completion of $C_0^\infty(\R^3)$ with respect to the metric of  $|x|\n \rot u\in L^2(\R^3)$ and $|x|\n\n\cdot u\in L^2(\mathbb{R}^3)$. \par We now introduce some preliminary results. \par
\begin{lemma}\label{PR}{\sl Let $D$ be an open and bounded subset of $\mathbb{R}^3$ having a smooth boundary $\partial D$. Let $u(t, x)\in C([0,T];W^{1,2}(D))$ with $u_t(t, x),D^2u(t, x)\in L^2(0,T;L^2(D))$. Let $g(x)$ be a smooth and bounded function with null trace on $\partial D$, then 
\be\label{PR-I}\intll 0t(u_\tau,\rot(g^2\rot u))d\tau= \frac{1}{2}\dm g\hskip0.02cm\rot u(t)\dm_2^2-\frac{1}{2}\dm g\hskip0.02cm\rot u(0)\dm_2^2\,,\mbox{ for all }t\in[0,T]\,.\ee}\end{lemma}\bp In the case of $g\equiv 1$ and  $u\in C([0,T];W_0^{1,2}(D))$ the result is classical and well known. In our assumptions on $u$ and $g$ the result holds following {e.g.} the technique employed in \cite[Lemma\,1]{Pro}. So we omit further details.\ep
Throughout the subsequent results, the function $\zeta$ will denote the function introduced in \eqref{eq:zeta}.
\begin{prop} \label{prop:relazione rotore-gradiente}
    {\sl Let $u\in W^{2,2}(\OO)\cap \mathcal V(\OO)$ and  $b:=\rot u$. Then, there exists a constant $c>0$, independent of $u$, such that
\be\label{RRG-I}
    \dm{\zeta\abs{x}\n b}\dm_2^2\le \dm{\zeta\abs{x}\rot b\dm_2^2+c\dm u}\dm_{1,2}^2,\ee
    for all $R>R_0$.}
\end{prop}\bp  The proof of \rf{RRG-I} follows from estimate \rf{EDM}. It is enough to set $h(x):=\zeta|x|$ and to take into account the properties of $\zeta$ and formula $$\n h=\n\zeta|x|+\zeta\n |x|\,,$$ for all $x\in \R^3$.\ep 
\begin{prop} \label{prop:J1}{\sl
    Let $u\in L^{\infty}(\Omega)$\, and $b\in W^{1,2}(\Omega)$. Then, given $\vep>0$,
    \begin{equation}
        |\big(u\cdot \n b,\zeta^2 |x|^2b\big)| \le \vep \dm{\zeta|x| \n b}\dm_2^2 + C(\vep)\dm{u}\dm_{\infty}^2\dm{\zeta|x| b}\dm_2^2\,.
    \end{equation}}
\end{prop}
\bp
The proof is an easy consequence of H\"older and Young inequalities.
\ep
\begin{prop} \label{prop:J2}{\sl}
  Let $u\in W^{1,2}(\Omega)\cap L^\infty(\Omega)$ and $b\in W^{1,2}(\Omega)$\,. Then, given $\vep>0$\,,
  \begin{equation}
     |\big(b\cdot\n u,\zeta^2|x|^2 b\big)|\le \vep\dm\zeta |x|\n b\dm_2^2+C(\vep)(\norm{b}_{2}^2+\dm u\dm_\infty^2\dm\zeta|x| b\dm_2^2)\,.
  \end{equation}
\end{prop}
\bp
Integrating by parts, we have
\[
\big(b\cdot\n u,\zeta^2|x|^2 b\big)=-\int_\Omega \left[(b\cdot \n (|x|^2\zeta^2)(b\cdot u) +2\zeta^2|x|^2(b\cdot\n b)\cdot u\right] \, dx\,.
\]
Hence, we should estimate
\[
\begin{aligned}
&H_1:= -2\int_{\Omega}\zeta (b\cdot u)|x|^2b\cdot \n \zeta \, dx\,,\\ 
&H_2:=-2\int_\Omega \zeta^2(x\cdot b)(b\cdot u)\,dx\,,\\
&H_3:=-2\int_\Omega \zeta^2 |x|^2(b\cdot \n b)\cdot u\, dx\,.
\end{aligned}
\]
Using H\"{o}lder and Young inequalities, considering that, for all $x\in\Omega$, $|\n \zeta(x)|\le c|x|^{-1}$\,, we get
\[
\ba{ll}
|H_1|\hskip-0.2cm&\le 2\norm{u}_{\infty}\dm{\zeta|x| b}\dm_2 \dm{\abs{\n \zeta}|x| b}\dm_2 \le c \dm{u}\dm_{\infty}^2\dm{\zeta |x|b}\dm_2^2+C\dm{b}\dm_{2}^2\,, \VS
|H_2|\hskip-0.2cm&\le2\dm u\dm_\infty \dm\zeta |x|b\dm_2  \dm\zeta b\dm_2\le c\dm u\dm_\infty^2 \dm\zeta |x|b\dm_2^2 +C\dm b\dm_{2}^2\,,\VS
    |H_3|\hskip-0.2cm&\le 2\dm{u}\dm_{\infty}\dm{\zeta|x| b}\dm_2 \dm{\zeta |x|\n b}\dm_2 \le \vep\dm{\zeta|x| \n b}\dm_2^2+  C(\vep) \norm{u}_{\infty}^2\dm{\zeta|x| b}\dm_2^2\,.
\ea
\]
Collecting the estimates, we prove the proposition.
\ep
\begin{prop} \label{prop:J3}
 Let $V=\xi+\omega\times x$, $(\xi,\omega)\in \mathbb{E}_3\times \mathbb{E}_3$, and $ b\in W^{1,2}(\Omega)$. Then, 
 \begin{equation}
     |(V\cdot\n b, \zeta^2\abs{x}^2b)| \le c\abs{\xi}^2\norm{\zeta \abs{x}b}_2^2+C\norm{b}_{2}^2\,,
 \end{equation}
 with $C,c>0$ independent of $V$ and $b$.
\end{prop}
\bp
Recalling that $\nabla\cdot V=0$, integrating by parts, we have
\[
\ba{ll}
\big(V\!\cdot\!\n b, \zeta^2 \abs{x}^2 b\big)\hskip-0.3cm&=\sfrac{1}{2} \big((\xi\!+\omega\times x), \zeta^2 \abs{x}^2 \n \abs{b}^2\big)\!=-\sfrac{1}{2} \big((\xi\!+\omega\times x), \!\n(\zeta^2 \abs{x}^2 ) \abs{b}^2\big)\\ \VSE=-\sfrac{1}{2} \big(\xi, \n(\zeta^2\abs{x}^2 ) \abs{b}^2\big)\,.
\ea
\]
In fact, by an elementary property of triple product, we have
\[
\big(\omega\times x, \n(\zeta^2 \abs{x}^2) \abs{b}^2\big)=0\,.
\]  Then,
\[
\begin{aligned}
& \sfrac{1}{2}\big(\xi, \n(\zeta^2 \abs{x}^2) \abs{b}^2\big)= \big(\xi,  \abs{x}^2 \zeta \cdot \n \zeta \abs{b}^2\big)+ \big(\xi, \zeta^2  \abs{b}^2x\big)\,.
\end{aligned}
\]
Hence, we deduce 
\[
\begin{aligned}
\sfrac{1}{2}\Big|\big(\xi,\! \n(\zeta^2 |x|^2) |b|^2\big)\Big|&\!\le  c |\xi|\dm \zeta|x|b\dm_2[\dm\abs{\n \zeta}\abs{x}b\dm_2\!+\dm{\zeta b}\dm_2]\!\VSE\!\le c\abs{\xi}^2\dm{\zeta\abs{x}b}\dm_2^2+C\dm b\dm_{2}^2 \,.
\end{aligned}
\]
The assertion is proved.
\ep
\begin{lemma} \label{le:Extra proprietà}
Let $(u,\pi,V)$ be a solution to \eqref{eq:model} corresponding to the initial datum $(u_0,\xi_0,\omega_0)\in \mathcal{V}(\Omega)\times\mathbb E_3\times \mathbb E_3$, satisfying \eqref{GSS-I}$_{1,3,4}$, \eqref{GSS-II} and \eqref{GSS-old}. Assume that
\[
\n u_0\in L^2(\Omega, \abs{x}).
\]
Then,
\begin{equation}
  \sup_{(0,T)}\dm  \abs{x}{\n u(t)} \dm_2+\intll0T\dm \abs{x}D^2u(t)\dm_2^2dt<\infty\,.
\end{equation}
\end{lemma}
\bp
We set $b(t,x):=\rot u(t,x)$. We multiply \eqref{eq:model}$_1$ by $A(t,x):=\rot\left[\zeta^2(x)\abs{x}^2 b(t,x)\right]$. Then,  integrating on $(0,T)\times\OO$, we get 
$$\intll0t(u_\tau(\tau,x),A(\tau,x))\,d\tau=\sfrac{1}{2}\dm 
\zeta(x) \abs{x}b(t)\dm_2^2-\sfrac{1}{2} \dm 
\zeta(x) \abs{x}b(0)\dm_2^2\,,\mbox{ for all }t\in[0,t)\,,$$
Where we employed  \rf{PR-I}, that we apply  by taking into account that $\zeta$ has
compact support, and  for $u$ \rf{GSS-I}$_{1,3,4}$-\rf{GSS-II} hold.
We consider the term
$$\intll0t(-\Delta u(\tau,x),A(\tau,x))\,d\tau=\intll0t(\zeta|x|\rot b(\tau),\zeta|x|\rot b(\tau))\,d\tau+\intll0t(\n(\zeta|x|)^2\times b(\tau),\rot b(\tau))\,d\tau\,.$$ Applying H\"order's inequality for the last integral, we get
$$\intll0t(-\Delta u(\tau,x),A(\tau,x))\,d\tau\geq \sfrac{1}{2}\intll0t\dm \zeta|x|\rot b(\tau)\dm_2^2\,d\tau-c\intll0t\dm u(\tau)\dm_{1,2}^2\,d\tau\,.$$
Recalling \rf{RRG-I}, we arrive at
\be\label{LP}\intll0t(-\Delta u(\tau,x),A(\tau,x))d\tau\geq \sfrac{1}{2}\intll0t\dm \zeta|x|\n b(\tau)\dm_2^2\,d\tau-c\intll0t\dm u(\tau)\dm_{1,2}^2\,d\tau\,.\ee
Furthermore, by definition, $A(t)\in J^2(\Omega_R)$ for a.a. $t\in (0,T)$, recalling  \rf{GSS-old} and the properties of $u$, then
\be\label{PT}
\intll0t(\nabla\pi(\tau), \rot(\zeta^2 \abs{x}^2b(\tau)))d\tau=0\,.
\ee
Finally, integration by parts yields to
\[
-\intll0t ((u-V)(\tau)\cdot \nabla u(\tau), \rot(\zeta^2 \abs{x}^2b(\tau)))\,d\tau\,=\,-\intll0t(\rot((u-V)(\tau)\cdot \nabla u(\tau)), \zeta^2 \abs{x}^2b(\tau))\,d\tau\,,
\]
and
\[
-\intll0t(\omega(\tau)\times u(\tau), \rot(\zeta^2 \abs{x}^2b(\tau)))\,d\tau\,=\,-\intll0t(\rot (\omega(\tau)\times u(\tau)), \zeta^2 \abs{x}^2b(\tau))\,d\tau\,.
\]
Hence, we have
\[
\ba {l} \displ
\overline{k}(t,R)+\int_0^t\overline{K}(\tau,R)\,d\tau \le \dm |x|\n u_0\dm_2^2+c(R_0)\int_0^t \norm{u(\tau)}_{2,2}^2\, d\tau \VS \hskip 2.2cm -2\int_0^t \Big[(\rot \left[(u-V)(\tau)\cdot \nabla u(\tau)\right], \zeta^2 \abs{x}^2b(\tau)) -(\rot(\omega(\tau)\times u(\tau)), \zeta^2\abs{x}^2b(\tau))\Big]\,d\tau,
\ea
\]
where
\[
\overline{k}(t,R):=\dm\zeta  \abs{x}b(t)\dm_{L^2(\OO^{ 2R_0})}^2,\mbox{ for all }t\in[0,T)\,,
\]
and
\[
\overline{K}(t,R):=\dm{\zeta \abs{x}\n b(t)}\dm_{L^2(\OO^{ 2R_0})}^2,\mbox{ a.e. in }t\in(0,T)\,,
\] where the dependence on $R$ is due to $\zeta$\,.
Now we estimate
\[
\int_0^t (\rot((u-V)(\tau)\cdot \nabla u(\tau)), \zeta^2\abs{x}^2 b(\tau) )\,d\tau \]
and \[\int_0^t (\rot(\omega(\tau)\times u(\tau)), \zeta^2 \abs{x}^2b(\tau))\,d\tau\,.
\]
We remark that
\[
\rot(u\cdot\n  u)=u\cdot \n b -b\cdot\n u, \quad \rot(V\cdot\n u)=V\cdot\n b-\omega\times b-\omega\cdot\n u, \quad \rot(\omega\times u)=-\omega\cdot\n u,
\]
then we should estimate
\[
\begin{aligned}
    &J_1:=\int_0^t \left(u(\tau)\cdot \n b(\tau),\zeta^2\abs{x}^2 b(\tau) \right)\,d\tau\,, \\
    &J_2:=\int_0^t \left(-b(\tau)\cdot\n u(\tau),\zeta^2 \abs{x}^2b(\tau)\right)\,d\tau\,,\\ 
    & J_3:= \int_0^t \left(-V(\tau)\cdot\n b(\tau), \zeta^2\abs{x}^2 b(\tau)\right)\,d\tau\,, \\
    &J_4:= \int_0^t \left(\omega(\tau)\times b(\tau), \zeta^2 \abs{x}^2b(\tau) \right)\,d\tau\,.
\end{aligned}
\]
By elementary arguments, we have
\[
J_4=0\,.
\]
The other integrals can be estimated easily by Propositions\,\ref{prop:J1}-\ref{prop:J3} and by energy inequality, we get
\[
\ba {rl}
\displ
\abs{J_1}\hskip-0.2cm&\displ\le \int_0^t \Big[ \vep \overline{K}(\tau,R)+ C\norm{u(\tau)}_{\infty}^2\overline{k}(\tau,R)+ c(R_0)\dm  u(\tau)\dm_{2,2}^2+c(R_0)\dm u(\tau)\dm_\infty^2\dm u(\tau)\dm_{1,2}^2\Big]\, d\tau\,,\VS
\abs{J_2}\hskip-0.2cm&\displ\le \!\int_0^t \Big[\vep\overline{K}(\tau,R)+C \norm{u(\tau)}_\infty^2\overline{k}(\tau,R)+c(R_0) \norm{u(\tau)}_\infty^2\norm{u(\tau)}_{1,2}^2+ c(R_0)\dm u(\tau)\dm_{2,2}^2\Big]\,d\tau\,,\VS
\abs{J_3}\hskip-0.2cm&\displ\le\int_0^t \Big[ca_0\overline{k}(\tau,R)+c(R_0)(a_0+1)\dm u(\tau)\dm_{1,2}^2\Big]\,d\tau\,, 
\ea
\]
where $a_0:=\norm{u_0}_2^2+\abs{\xi_0}^2+\abs{\omega_0}^2$\,.
Collecting all the estimates, choosing $\vep>0$ sufficiently small and using Proposition\,\ref{prop:relazione rotore-gradiente}, we find that
\[
\begin{aligned}
&\overline{k}(t,R)+\sfrac{1}{2}\int_0^t \overline{K}(\tau,R)\, d\tau \le \dm |x|\n u_0\dm_2^2\\ &\hskip1cm+C(R_0,a_0)\int_0^t\Big[ (1+ \norm{u(\tau)}_\infty^2) \overline{k}(\tau,R)+\dm u(\tau)\dm_\infty^2\norm{ u(\tau)}_{1,2}^2+\dm u\dm_{2,2}^2\Big]\, d\tau \,.
\end{aligned}
\]
Hence, employing  Gronwall lemma, letting   $R\to\infty$, by virtue of Beppo-Levi theorem, we arrive at
$$\dm |x|\n u(t)\dm_{L^2(\OO^{2R_0})}^2+\intll0t\dm |x|\n\rot u(\tau)\dm_{L^2(\OO^{2R_0})}^2\,d\tau\leq c(u_0,\xi_0,\omega_0,t)\,.$$ Since by the assumptions on $u$ the contribute on $\OO\cap B_{2R_0}$ of the previous norms is finite, we have proved 
\[
\sup_{t\in(0,T)}\norm{\abs{x}\rot u(t)}_2^2+\int_0^T \norm{\abs{x}\n\rot u(t)}_2^2\, dt <\infty.
\]
By estimates \rf{RGW} and \rf{RDSW}, we can conclude the proof.
\ep
\begin{lemma} \label{lemmma:Galdipreliminari}
    Let $u\in\mathcal{V}(\Omega)$. Then,
    \begin{equation} \label{eq:gradsimmetrico}
        \norm{\n u}_2=\sqrt{2}\norm{\mathbb{D}(u)}_2.
    \end{equation}
    Moreover, there exists $C>0$, independent of $u$,  such that
    \begin{equation} \label{eq:CM}
    \norm{u}_3\le C\norm{u}_2^{\frac{1}{2}} \norm{\n u}_2^{\frac{1}{2}}.
    \end{equation}
\end{lemma}
\bp
For the proof of \eqref{eq:gradsimmetrico} we quote \cite{Ga:Hand}, for the proof of \eqref{eq:CM} we quote \cite{CM}, where it is presented in a more general setting.
\ep
\section{Proof of the main results}
\bp[Proof of Theorem\,\ref{TDU}]
Property \eqref{TDU-I} follows immediately from Lemma\,\ref{le:Extra proprietà}. \par
We now show that also \eqref{TDU-II} holds. \par
Multiply \eqref{eq:model}$_1$ by $\varphi_R^2u_t$, where $\varphi_R$ is a usual cut-off function and $R>2\bar{R}:=2\text{diam}(\mathcal{B}) $. We get
\begin{equation*}
\ba {l} \displ
 \sfrac{d}{dt}\dm\varphi_R\mathbb{D}(u)\dm_{2}^2+|{\dot{\xi}}|^2+\dot{\omega}\cdot I\cdot \dot{\omega}+\norm{\varphi_Ru_t}_{2}^2=-\omega\times\xi\cdot \dot{\xi}-\omega\times(I\cdot \omega)\cdot\dot{\omega}-(u\cdot\n u, \varphi_R^2 u_t)\VS \hskip7cm +(V\cdot\n u,\varphi_R^2 u_t)-(\omega\times u, \varphi_R^2 u_t) - (\mathbb{T}\cdot \varphi_R\n \varphi_R, u_t)\,.
\ea
\end{equation*}
We now estimate the terms on the right-hand side of the previous equation. \par By H\"{o}lder and Young inequalities, lemma\,\ref{lemmma:Galdipreliminari} and energy estimates, we have
\[
\begin{array}{rl}
|(u\cdot\n u, \varphi_R^2u_t)| \hskip-0.2cm&\le \norm{u}_\infty\norm{\varphi_R\n u}_2\norm{\varphi_Ru_t}_2\le \vep  \norm{\varphi_Ru_t}_2^2+C(\vep)\norm{\varphi_R\mathbb{D}( u)}_2^2\norm{u}_{\infty}^2\,, \VS
|(V\cdot\n u, u_t)|\hskip-0.2cm&=|((\xi+\omega\times x)\cdot \n u,\varphi_R^2 u_t)|\VSE\le \abs{\xi}\norm{\varphi_R\n u}_2\norm{\varphi_Ru_t}_2+ \abs{\omega}\norm{\varphi_R\abs{x}\abs{\n u}}_2\norm{\varphi_Ru_t}_2
\VSE\leq \vep \norm{\varphi_Ru_t}_2^2 + C(\vep)(\dm u_0 \dm_2^2+\abs{\xi_0}^2+\abs{\omega_0}^2)(\norm{\varphi_R\abs{x}\abs{\n u}}_2^2+\norm{\varphi_R\mathbb{D} (u)}_2^2)\,,\VS
|(\omega\times u, \varphi_R^2u_t)|\hskip-0.2cm&\le c\abs{\omega}\norm{\varphi_Ru}_2\norm{\varphi_Ru_t}_2 \le \vep \norm{\varphi_Ru_t}_2^2+ C(\vep)(\abs{\xi_0}^2+\abs{\omega_0}^2+\norm{u_0}_2^2)^2\,,\VS
|{\omega\times\xi\cdot \dot{\xi} }|\hskip-0.2cm&\le \vep |\dot{\xi}|^2+C(\vep)(\abs{\xi_0}^2+\abs{\omega_0}^2+\norm{u_0}_2^2)^2\,,\VS
\abs{\omega\times(I\cdot \omega)\cdot\dot{\omega}} \hskip-0.2cm&\le \vep \abs{\dot{\omega}}^2+C(\vep)(\abs{\xi_0}^2+\abs{\omega_0}^2+\norm{u_0}_2^2)^2\,.
\end{array}
\]
Finally, we estimate
\[
\abs{\int_\Omega \varphi_R\n \varphi_R \cdot \mathbb{T}\cdot u_t} \le c \int_{R\le \abs{x}\le 2R}\varphi_R\abs{x}^{-1}\abs{\mathbb{T}}\abs{u_t} \le c \int_{R\le \abs{x}\le 2R}\abs{x}^{-1}(\abs{\pi}+\abs{\mathbb{D}(u)})\varphi_R\abs{u_t}.
\]
By Young inequality, we have
\[
\begin{aligned}
& \int_{R\le \abs{x}\le 2R}\abs{x}^{-1}(\abs{\pi}+\abs{\mathbb{D}(u)})\varphi_R\abs{u_t}\le \vep\norm{\varphi_Ru_t}_2^2+ C(\vep)\int_{\Omega^R}\abs{x}^{-2}(\abs{\pi}^2+\abs{\mathbb{D}(u)}^2)=: \vep\norm{\varphi_Ru_t}_2^2+I_1(R)\,.
\end{aligned}
\]
 Hence, choosing $\vep>0$ sufficiently small, we get the following inequality
\[
\ba{l}
\sfrac{d}{dt}\norm{\varphi_R\mathbb{D}(u)}_2^2\!+\frac{1}{2}\big[|\dot{\xi}|^2\!+\dot{\omega}\cdot\! I\!\cdot \dot{\omega}+\norm{\varphi_Ru_t}_2^2\big] \le C(\vep)\big[\abs{\xi_0}^2\!+\abs{\omega_0}^2\!+\norm{u_0}_2^2\!+\norm{u}_\infty^2\!+\norm{\varphi_R u}_2^2\big]\norm{\varphi_R\mathbb{D}(u)}_2^2\VS \hskip 3.8cm +C(\vep)\big[\abs{\xi_0}^2\!+\abs{\omega_0}^2\!+\norm{u_0}_2^2\big]\norm{\varphi_R\abs{x}\n u}_2^2\!+C(\vep)\big[\abs{\xi_0}^2\!+\abs{\omega_0}^2\!+\norm{u_0}_2^2\big]^2\!+I_1(R).
\ea
\]
By the integrability properties of $\mathbb{D}(u)$, the multidimensional Hardy inequality and Lebesgue theorem, we can say that
\[
\lim_{R\to+\infty}I_1(R):=\lim_{R\to+\infty}\int_{\Omega^R}\abs{x}^{-2}(\abs{\pi}^2+\abs{\mathbb{D}(u)}^2)=0.
\]
Therefore, integrating on $(0,t)$ and letting $R\to+\infty$, by virtue of Beppo-Levi theorem, we deduce
\[
\ba{l} 
\displ
\norm{\mathbb{D}(u(t))}_2^2+\sfrac{1}{2}\intll0t\Big[|\dot{\xi}|^2+\dot{\omega}\cdot I\cdot \dot{\omega}+\norm{u_t}_2^2\Big]d\tau \le \norm{\mathbb{D}(u(0))}_2^2\VS \hskip 3cm +C(\vep)(\abs{\xi_0}^2\!+\abs{\omega_0}^2\!+\norm{u_0}_2^2)\intll0t\norm{\mathbb{D}(u)}_2^2d\tau+C(\vep)\intll0t \norm{u}_{\infty}^2\norm{\mathbb{D}(u)}_2^2d\tau\VS \hskip 4cm +C(\vep)(\abs{\xi_0}^2\!+\abs{\omega_0}^2\!+\norm{u_0}_2^2)\intll0t\norm{\abs{x}\n u}_2^2d\tau+C(\vep)(\abs{\xi_0}^2\!+\abs{\omega_0}^2\!+\norm{u_0}_2^2)^2t .
\ea
\]
This completes the proof.
\ep
\bp[Proof of Theorem\,\ref{UT}]
Assume that $(u,\pi,V)$ and $(\bar{u},\bar{\pi},\bar{V})$ are two solutions corresponding to the same initial datum $(u_0,\xi_0,\omega_0)$, where $u$ enjoys the property \eqref{UT-I} and $(\bar{u},\bar{\pi},\bar{V})$ satisfies the energy inequality. Set
\[
\widehat {u}:=u-\bar{u}, \quad \widehat {\pi}:= \pi-\bar{\pi}, \quad\widehat {V}:=V-\bar{V}, \quad \widehat {\xi}:=\xi-\bar{\xi}, \quad \widehat {\omega}:=\omega-\bar{\omega}.
\]
Then, $(\widehat {u},\widehat {\pi},\widehat {V})$ satisfies
\begin{equation} \label{eq:differenza}
\begin{cases}
    \widehat {u}_t=\nabla\cdot\mathbb{T}(\widehat {u},\widehat {\pi})-(\widehat {u}-\widehat {V})\cdot\n u-(u-V)\cdot \n \widehat {u} +(\widehat {u}-\widehat {V})\cdot \n \widehat {u}-\widehat {\omega}\times u-\bar{\omega}\times \widehat {u}\,,\\ 
\nabla\cdot \widehat {u}=0\,,\quad \forall(t,x)\in (0,T^*)\times\Omega\,,\\
\widehat {u}(t,x)=\widehat {\xi}(t)+\widehat {\omega}(t)\times x\,, \quad \forall(t,x)\in (0,T^*)\times\partial\Omega\,,\\
 \lim_{\abs{x}\to \infty} \widehat {u}(t,x)=0 \,,\quad \forall t\in (0,T)\,, \\
        \dot{\widehat {\xi}} + \widehat {\omega} \times \widehat {\xi} +  \int_{\partial\Omega} \mathbb{T}(\widehat {u}, \widehat {\pi})\cdot \nu=0 \,,\quad \forall t\in (0,T)\,,\\
         I \cdot \dot{\widehat {\omega}}+ \widehat {\omega} \times(I\cdot \widehat {\omega})+ \int_{\partial\Omega} x \times \mathbb{T}(\widehat {u},\widehat {\pi}) \cdot \nu=0 \,,\quad \forall t\in (0,T)\,,\\
         \widehat {\xi}(0)= 0\,, \quad \widehat {\omega}(0)=0\,, \\
         \widehat {u}(0,x)=0\,, \quad \forall x \in \Omega\,.
    \end{cases}
\end{equation}
Multiply \eqref{eq:differenza}$_1$ by $\widehat {u}$ and integrate over $\Omega$. Since 
\[
\ba{rl}
(\bar{\omega}\times \widehat {u},\widehat {u})\hskip-0.2cm& =0\,,\VS
((u-V)\cdot\n \widehat{u},\widehat{u}) \hskip-0.2cm&=0\,,\VS
((\widehat{u}-\widehat{V})\cdot\n \widehat {u},\widehat {u})\hskip-0.2cm &=0\,,
\ea
\] 
we get
\[
\begin{aligned}
&\sfrac{1}{2}\sfrac{d}{dt}(\norm{\widehat {u}}_2^2+|{{\widehat {\xi}}}|^2+\widehat {\omega}\cdot I\cdot \widehat {\omega})+\norm{\mathbb{D}(\widehat {u})}_2^2=\widehat {\xi}\times\widehat {\omega}+(I\cdot\widehat {\omega})\times\widehat {\omega}-((\widehat {u}-\widehat {V})\cdot\n u,\widehat {u})-(\widehat {\omega}\times u,\widehat {u}).
\end{aligned}
\]
We estimate the right-hand side terms using H\"{o}lder and Young inequalities. We find
\[
\ba {rl}
|{(\widehat {\omega}\times u,\widehat {u})}|\hskip-0.2cm &\le c|{\widehat {\omega}}|\dm{u}\dm_2\dm{\widehat {u}}\dm_2\le c|{\widehat {\omega}}|^2+c\dm{u}\dm_2^2\dm{\widehat {u}}\dm_2^2\,, \VS
|{\widehat  {\xi}\times\widehat {\omega}}|\hskip-0.2cm&\le c|{\widehat {\xi}}|^2+c|{\widehat {\omega}}|^2\,, \VS
|{\widehat {\omega}\cdot I\cdot \widehat {\omega}}|\hskip-0.2cm&\le c|{\widehat {\omega}}|^2\,,\VS
|\widehat {\xi}\times\widehat {\omega}+(I\cdot\widehat {\omega})\times\widehat {\omega}|\hskip-0.2cm&\le c|\widehat {\xi}|^2+c|\widehat {\omega}|^2\,, \VS
((\widehat {u}-\widehat {V})\cdot\n u,\widehat {u}) \hskip-0.2cm&=(\widehat {u}\cdot\n u,\widehat {u})-(\widehat {V}\cdot\n u,\widehat {u})=:A_1-A_2\,.
\ea
\]
We can easily estimate $A_1$ and $A_2$ by H\"{o}lder and Young inequalities and inequality \eqref{eq:CM}, in fact
\[
|A_1| \le \norm{\widehat {u}}_6 \norm{\n u}_2\norm{ \widehat {u}}_3 \le \vep\norm{\n \widehat {u}}_2^2+ C(\vep)\norm{\n u}_2^4 \norm{\widehat {u}}_2^2\,,
\]
\[
|A_2|\le |{\widehat {\xi}}|\norm{\n u}_2\norm{\widehat {u}}_2+|{\widehat {\omega}}|\norm{\abs{x}\abs{\n u}}_2\norm{\widehat {u}}_2\le C(|{\widehat {\xi}}|^2+\abs{\widehat {\omega}}^2)+C(\norm{\n u}_2^2+\norm{\abs{x}\abs{\n u}}_2^2)\norm{\widehat {u}}_2^2\,.
\]
Collecting all the estimates and choosing $\vep>0$ sufficiently small, setting
\[
y(t):=\norm{\widehat {u}}_2^2+|{{\widehat {\xi}}}|^2+\widehat {\omega}\cdot I\cdot \widehat {\omega},
\]
we find
\[
\sfrac{d}{dt}y(t)\le C\varphi(t)y(t),
\]
with $\varphi\in L^1(0,T)$. Hence, by virtue of Gronwall lemma, function $y(t)$ is  null, that is the uniqueness of $(u,\pi,V)$\,.
\ep

{\bf Acknowledgements.}
 \small{The authors acknowledge the support from GNFM research group of the \textit{Istituto Nazionale di Alta Matematica}.}

\end{document}